\documentclass[10pt]{IEEEtran}

\ifCLASSINFOpdf
\else
\fi
\usepackage{amssymb}
\usepackage{multirow}
\usepackage{indentfirst}
\usepackage{inputenc}
\usepackage{array}
\usepackage[switch,mathlines]{lineno}
\usepackage{soul}
\usepackage{graphicx}
\usepackage{cite}
\usepackage[cmex10]{amsmath}
\usepackage{algorithmic}
\usepackage{array}
\usepackage{fixltx2e}
\usepackage{CJK}
\usepackage{algorithmic}
\usepackage{array}
\usepackage{fixltx2e}
\usepackage{bm}
\usepackage{color}
\usepackage{amsmath,amsthm,amssymb}
\newtheorem{theorem}{Theorem}
\newtheorem{corollary}{Corollary}[theorem]

\newtheorem*{remark}{Remark}
\begin{document}
%
\title{Robust Power System Dynamic State Estimator with Non-Gaussian Measurement Noise: Part I--Theory}
\author{Junbo~Zhao,~\IEEEmembership{Student Member,~IEEE},~Lamine~Mili,~\IEEEmembership{Fellow,~IEEE}\\

\thanks{Junbo Zhao and Lamine Mili are with the Bradley Department of Electrical Computer Engineering,
Virginia Polytechnic Institute and State University, Northern Virginia Center, Falls Church, VA 22043, USA (e-mail: zjunbo@vt.edu, lmili@vt.edu).}}
\markboth{IEEE TRANSACTIONS ON POWER SYSTEMS,~Vol.~, No.~, ~2017}%
{Shell \MakeLowercase{\textit{et al.}}: Bare Demo of IEEEtran.cls for Journals}
\maketitle

\begin{abstract}
This paper develops the theoretical framework and the equations of a new robust Generalized Maximum-likelihood-type Unscented Kalman Filter (GM-UKF) that is able to suppress observation and innovation outliers while filtering out non-Gaussian measurement noise. Because the errors of the real and reactive power measurements calculated using Phasor Measurement Units (PMUs) follow long-tailed probability distributions, the conventional UKF provides strongly biased state estimates since it relies on the weighted least squares estimator. By contrast, the state estimates and residuals of our GM-UKF are proved to be roughly Gaussian, allowing the sigma points to reliably approximate the mean and the covariance matrices of the predicted and corrected state vectors. To develop our GM-UKF, we first derive a batch-mode regression form by processing the predictions and observations simultaneously, where the statistical linearization approach is used. We show that the set of equations so derived are equivalent to those of the unscented transformation. Then, a robust GM-estimator that minimizes a convex Huber cost function while using weights calculated via Projection Statistics (PS's) is proposed. The PS's are applied to a two-dimensional matrix that consists of serially correlated predicted state and innovation vectors to detect observation and innovation outliers. These outliers are suppressed by the GM-estimator using the iteratively reweighted least squares algorithm. Finally, the asymptotic error covariance matrix of the GM-UKF state estimates is derived from the total influence function. In the companion paper, extensive simulation results will be shown to verify the effectiveness and robustness of the proposed method.
\end{abstract}
\vspace{-0.2cm}

%
\IEEEpeerreviewmaketitle
\vspace{-0.5cm}
\section{Introduction}
\IEEEPARstart{T}{he} widespread deployment of synchro-Phasor Measurement Units (PMUs) on power transmission grids has made possible the real-time monitoring and control of power system dynamics. However, these functions cannot be reliably achieved without the development of a fast and robust Dynamic State Estimator (DSE). Indeed, the state variable estimates of the synchronous machines can be utilized by power system stabilizers, automatic voltage regulators, and under-frequency relays to enhance small signal stability and to initiate generation outages and load shedding during transient instabilities, among other actions \cite{Kamwa2001,Mili2002}.

\begin{figure*}[!t]
\centering
\includegraphics[height=16cm]{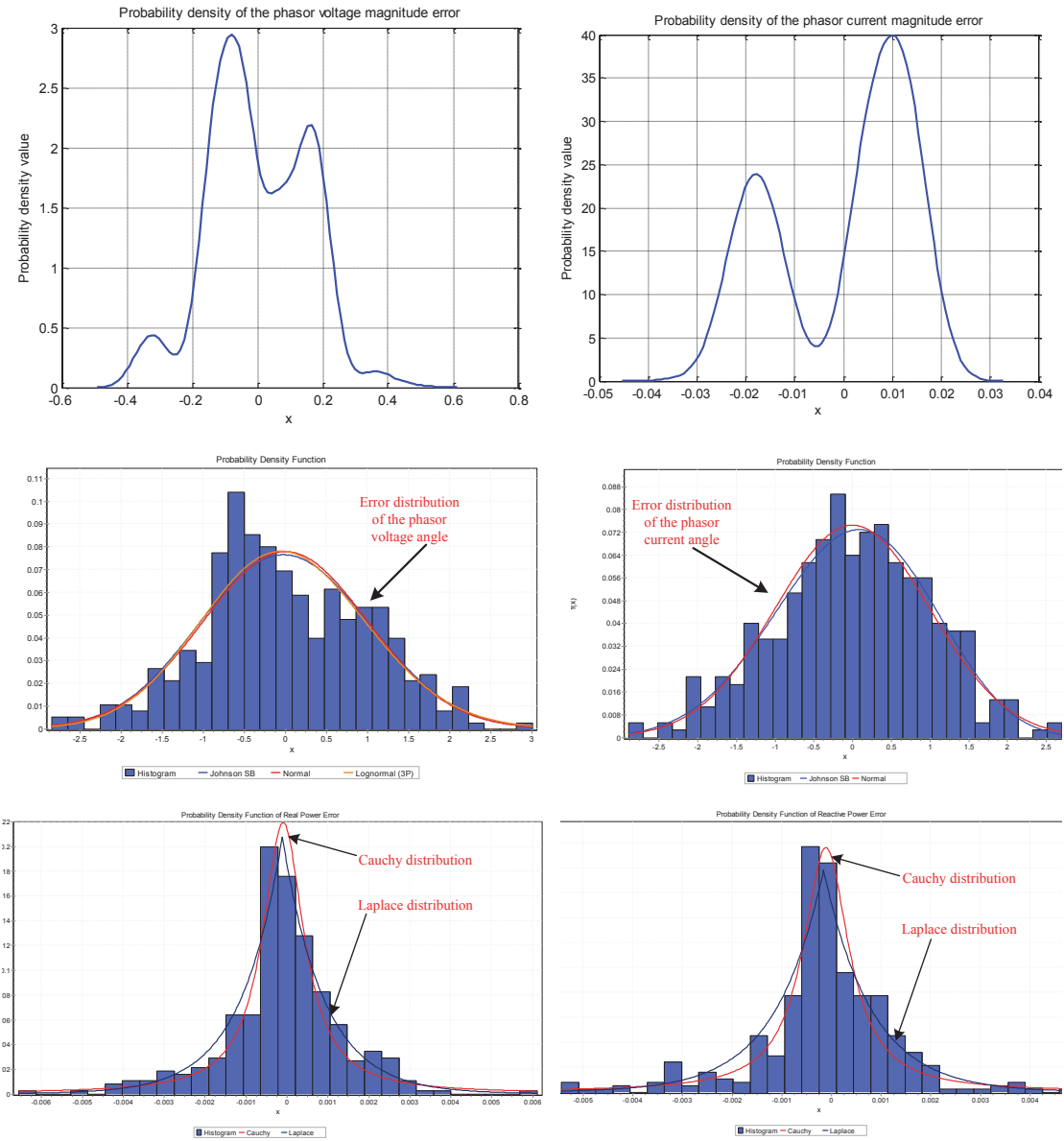}
\caption{Error distributions of the phasor voltage angle and magnitude, the phasor current angle and magnitude, the real and reactive power using field PMU data.}
\label{Fig.non_gaussian_noise_real_data}
\end{figure*}
To date, a variety of dynamic state estimators have been proposed in the literature; they are based on the Extended Kalman Filter (EKF) \cite{EKFKamwa2011,ZHuangKJ2007}, the Iterated EKF (IEKF) \cite{FanWhbe2013,Junbo_GMIEKF2016}, the unscented Kalman filter (UKF) \cite{ABCPal2014,APSMelipoulos2012,Rouhani2016}, to cite a few. However, all these methods suffer from several important shortcomings, precluding them from being adopted by power utilities for power system real-time applications. To be specific, they cannot handle \emph{i) non-Gaussian process and observation noise of the system nonlinear dynamic models and ii) innovation, observation and structural outliers}.

There are several reasons for these shortcomings. Firstly, the current DSE approaches assume that both the process and the observation noise of the system nonlinear dynamic models are Gaussian. However, two recent investigations conducted by PNNL \cite{henry2015_PNNL,PNNL_report2014} revealed that the PMU measurement errors of the voltage and current magnitudes obey non-Gaussian probability distributions. This is demonstrated in Fig.\ref{Fig.non_gaussian_noise_real_data} using real PMU data provided to us by PNNL. This figure displays histograms and parametric probability density estimates of PMU errors on nodal voltage magnitudes and angles, line current magnitudes and angles, and line real and reactive powers. As observed in Fig.\ref{Fig.non_gaussian_noise_real_data}, except for the measurement errors on nodal voltage and line current angles, which are roughly Gaussian, the measurement errors on both nodal voltage and line current magnitudes obey a bimodal Gaussian mixture distribution. As for the measurement errors of line real and reactive powers calculated from voltage and current phasors, they follow a thick tailed distribution that may be approximated by either the Laplacian or the Cauchy distribution. Recall that in contrast to the Gaussian distribution, which is a short-tailed distribution, a thick-tailed distribution is the one that allows the associated random variable to take, as compared to a scale parameter, large values with a non-negligible probability. Evidently, the presence of non-Gaussian noise calls for new research and development in robust power system DSE based on robust statistics.

Secondly, three types of outliers associated with a given dynamical system model have been defined by Gandhi and Mili \cite{Lmili2010}, namely observation outliers, which affect the metered values; innovation outliers, which corrupt the predicted state estimates; and structural outliers, which affect the system dynamic states and the observation functions. Observation outliers may result from large biases in PMU measurements due to infrequent calibration, or instrument failures, or impulsive communication noise \cite{Martin2007,KAThorp2015}. As for innovation outliers, they may occur in several different ways. For example, some of the generator models may not be well calibrated, resulting in highly inaccurate model outputs that are inconsistent with the measurements. This was precisely the case in the 1996 blackout, where the model being used predicted system stability while in reality the system was undergoing numerous cascading failures, which resulted in a rapid system collapse that occurred within minutes \cite{Kosterev1999,Kosterev2013}. Innovation outliers may also be induced by the approximations in the state prediction model or by a system process impulsive noise. By contrast, structural outliers are induced by wrong circuit breaker statuses or gross errors in the model parameters of the transmission lines, or of the automatic voltage regulators, or of the synchronous machines. In \cite{Pal2015_DSE}, it is reported that wrong estimates of the parameters of the synchronous machine models may result from the use of erroneous metered values. It turns out that the conventional filters, namely the EKF, the IEKF, the UKF, and the Particle Filter (PF) are not robust to any type of outliers. For instance, it is demonstrated in \cite{NZhou2015} that their performances are significantly degraded in the presence of observation outliers. To address this issue, Rouhani and Abur \cite{Rouhani2016} developed a robust UKF-based DSE using the Least-Absolute-Value (LAV) estimator. However, the authors do not address the vulnerability of the DSE to innovation outliers. In \cite{Junbo_GMIEKF2016}, a robust IEKF was proposed to handle observation and innovation outliers, but it may suffer from divergence problems if the nonlinearity of the system model is strong. In addition, both \cite{Rouhani2016,Junbo_GMIEKF2016} do not address the non-Gaussianity of the measurement noise.

In this paper, a robust Generalized Maximum-Likelihood-type UKF (GM-UKF) method is proposed to suppress observation and innovation outliers while filtering out non-Gaussian measurement noise. Our choice of the UKF is motivated by the fact that, considering the real-time implementation requirements for nonlinear DSE, it achieves a more balanced performance between computational efficiency and ability to cope with strong system nonlinearities than the EKF, or the IEKF, or the PF. However, the UKF is based on the sigma points, which reliably approximate the mean and the covariance matrices of the state estimates only under the Gaussian assumption of the process and observation noises. We show that this assumption is further stressed by the reliance of the UKF on the weighted least squares estimator. Interestingly, the state estimates calculated by our GM-UKF are shown to be asymptotically Gaussian even when the noises obey thick-tailed distributions, which is precisely the case when using PMU measurements. Furthermore, we show that the state estimates obtained from the application of statistical linearization to the nonlinear discrete-time state space system model are equivalent to those of the unscented transformation. Therefore, our filter allows the sigma points to provide good results.

It is developed according to the following steps. We first derive a redundancy batch-mode regression form by processing the predictions and observations simultaneously; this overdetermined system of equations provides the data redundancy needed for the detection and suppression of the innovation and observation outliers. This is achieved by means of a robust GM-estimator defined as the minimum of the Huber convex cost function while using weights calculated via the Projection Statistics (PS's). The latter are applied to a two-dimensional matrix consisting of serially correlated predicted state and innovation vectors. Then, a statistical test is applied to them to flag the outliers. Finally, the GM-estimator is solved via the iteratively reweighted least squares algorithm and the asymptotic error covariance matrix of the state estimates is calculated from the total influence function.

The rest of the paper is organized as follows. Section II presents the problem formulation. Section III develops the theory of the proposed GM-UKF and finally Section IV concludes the paper.
\vspace{-0.2cm}
\section{Problem Formulation}
\subsection{Nonlinear Discrete-Time Dynamical System Model}
A discrete-time state space representation of a general nonlinear dynamical system is expressed as
\begin{equation}
{\bm{x}_k} = \bm{f}\left( {{\bm{x}_{k - 1}},{{\bm{u}}_{k}}} \right) + {{\bm{w}}_{k}},
\label{Eq:discrete_state_model}
\end{equation}
\vspace{-0.3cm}
\begin{equation}
\bm{z}_k = {\bm{h}}\left( {{\bm{x}_k},{{\bm{u}}_k}} \right) + {{\bm{v}}_k},
\label{Eq:discrete_observation_model}
\end{equation}
where $\bm{x}_k\in \mathbb{R}^{n\times 1}$ and $\bm{z}_k\in \mathbb{R}^{m\times 1}$ are the state vector and the measurement/observation vector at time sample $k$, respectively; $\bm{f}$ and $\bm{h}$ are vector-valued nonlinear functions; $\bm{w}_k$ and $\bm{v}_{k}$ are the system process and observation noise, respectively; they are assumed to be independent and identically distributed with zero mean and covariance matrices $\bm{Q}_k$ and $\bm{R}_k$, respectively; $\bm{u}_k$ is the system input vector.
\vspace{-0.3cm}
\subsection{Dynamic State Estimation using UKF}
The main idea underlying the UKF is the application of a deterministic sampling technique known as the unscented transformation, which allows us, under the Gaussian noise assumption, to choose a set of sample points, termed sigma points, that have the same mean and covariance matrix as those of the a priori state vector \cite{Julier2000}. These sigma points are then propagated through the non-linear functions $\bm{f}$ and $\bm{h}$, yielding an estimation of the a posteriori state statistics by using the Kalman filter approach, i.e., the sample mean and the sample covariance matrix. Consequently, no calculation of Jacobian matrices is required, which can be by itself a difficult task to achieve in some cases or computationally costly.

To be specific, given a state estimate at time step $k$-1, ${\bm{\widehat {x}}_{k - 1\left| {k - 1} \right.}}\in \mathbb{R}^{n\times 1}$, having a covariance matrix given by ${\bm{P}_{k - 1|{k - 1}}^{xx}}$, its statistics are captured by 2$n$ weighted sigma points defined as
\begin{equation}
\bm{\chi} _{_{k - 1\left| {k - 1} \right.}}^i = {\bm{\widehat x}_{k - 1\left| {k - 1} \right.}} \pm {\left( {\sqrt {n{\bm{P}_{k - 1\left| {k - 1} \right.}^{xx}}} } \right)_i},
\label{Eq:sigma_points1}
\end{equation}
with weights ${w_i} = \frac{1}{{2n}},i=1,...,2n$. Then, each sigma point is propagated through the nonlinear system process model (\ref{Eq:discrete_state_model}), yielding a set of transformed samples expressed as
\begin{equation}
\bm{\chi}_{_{k\left| {k - 1} \right.}}^i = \bm{f}\left( {\bm{\chi}_{_{k - 1\left| {k - 1} \right.}}^i} \right).
\label{Eq:propagated_sigma_points}
\end{equation}
Next, the predicted sample mean and sample covariance matrix of the state vector are calculated by
\begin{equation}
{\bm{\widehat x}_{k\left| {k - 1} \right.}} = \sum\limits_{i = 1}^{2n} {{w_i}} \bm{\chi}_{_{k\left| {k - 1} \right.}}^i,
\label{Eq:propagated_sigma_points_mean}
\end{equation}
\vspace{-0.3cm}
\begin{equation}
{\bm{P}_{k| {k - 1}}^{xx}} = \sum\limits_{i = 1}^{2n} {{w_i}}( {\bm{\chi}_{_{k| {k - 1}}}^i - {{\bm{\widehat x}}_{k| {k - 1}}}}){( {\bm{\chi}_{_{k|{k - 1}}}^i - {{\bm{\widehat x}}_{k|{k - 1}}}})^T}+\bm{Q}_k.
\label{Eq:propagated_sigma_points_covariance}
\end{equation}
Finally, the measurement updating is performed and the filtered state ${\bm{\widehat x}_{k\left| k \right.}}$ with the covariance matrix ${\bm{P}_{k\left| k \right.}^{xx}}$ are calculated by
\begin{equation}
{\bm{K}_k} = \bm{P}_{_{k\left| {k - 1} \right.}}^{xz}{\left( {\bm{P}_{_{k\left| {k - 1} \right.}}^{zz}} \right)^{ - 1}},
\label{Eq:ukf_gain}
\end{equation}
\begin{equation}
{\bm{\widehat x}_{k\left| k \right.}} = {\bm{\widehat x}_{k\left| {k - 1} \right.}} + {\bm{K}_k}\left( {{\bm{z}_{k}} - {{\bm{\widehat z}}_{k\left| {k - 1} \right.}}} \right),
\label{Eq:ukf_filtering}
\end{equation}
\begin{equation}
{\bm{P}_{k\left| k \right.}^{xx}} = {\bm{P}_{k\left| {k - 1} \right.}^{xx}} - {\bm{K}_k}\bm{P}_{_{k\left| {k - 1} \right.}}^{zz}\bm{K}_k^T,
\label{Eq:ukf_covariance_matrix_updating}
\end{equation}
where ${\bm{\widehat z}_{k| {k - 1}}} = \sum\limits_{i = 1}^{2n} {{w_i}} \bm{z}_{_{k|{k - 1}}}^i$ is the predicted measurement vector and $\bm{z}_{_{k| {k - 1}}}^i = \bm{h}({\bm{\chi}_{_{k| {k - 1}}}^i})$; the self and cross-covariance matrices, $\bm{P}_{_{k| {k - 1}}}^{zz}$ and $\bm{P}_{_{k| {k - 1}}}^{xz}$, are respectively calculated by
\begin{equation}
\bm{P}_{_{k| {k - 1}}}^{zz} = \sum\limits_{i = 1}^{2n} {{w_i}}( {\bm{z}_{_{k| {k - 1}}}^i - {{\bm{\widehat z}}_{k| {k - 1}}}}){( {\bm{z}_{_{k| {k - 1}}}^i - {{\bm{\widehat z}}_{k|{k-1}}}})^T} + {\bm{R}_k},
\label{Eq:ukf_cross_covariance_matrix_zz}
\end{equation}
\begin{equation}
\bm{P}_{_{k| {k - 1}}}^{xz} = \sum\limits_{i = 1}^{2n} {{w_i}}( {\bm{\chi}_{_{k| {k - 1}}}^i - {{\bm{\widehat x}}_{k| {k - 1}}}}){( {\bm{z}_{_{k| {k - 1}}}^i - {{\bm{\widehat z}}_{k|{k-1}}}})^T}.
\label{Eq:ukf_cross_covariance_matrix_xz}
\end{equation}
\subsection{Motivation of the Use of a Robust UKF}
If the system process and measurement noises obey a Gaussian probability distribution, the filtered state, ${\bm{\widehat {x}}_{k - 1| {k -1}}}$, will follow a Gaussian distribution as well. In that case, the sample mean and the sample covariance matrix of ${\bm{\widehat {x}}_{k - 1| {k -1}}}$ will be captured by the sigma points and the UKF will produce reliable state estimates. However, the Gaussianity assumption may not hold true in practice. This is precisely the case in power systems; for instances, impulsive process noise may occur due to system model inaccuracy at a certain time window and the PMU measurement noise may not follow a Gaussian distribution as shown in Fig.\ref{Fig.non_gaussian_noise_real_data}. Consequently, the sigma points may not capture the complete statistics of the state vector, resulting in poor or even diverged estimations. Furthermore, since the UKF lacks statistical robustness, it is sensitive to any type of outliers, including observation, innovation and structural outliers. In power system DSE, observation outliers refer to the phase biases and gross errors in PMU measurements \cite{KAThorp2015}; innovation outliers may be induced by incorrect generator parameter values, failure of brushless exciter rotating diodes, or impulsive system process noise; and structural outliers may be caused by transmission parameter errors or topology errors. In the following section, we will propose a robust GM-UKF that is able to suppress observation and innovation outliers and to filter out various types of thick-tailed measurement noises. Note that the problem of the identification and suppression of structural outliers is outside the scope of this paper since it requires a different formulation; it will be addressed in a future work.
\vspace{-0.2cm}
\section{The Proposed GM-UKF}
Our GM-UKF consists of four major steps, namely a batch-mode regression form step, a robust pre-whitening step, a robust regression state estimation step, and a robust error covariance matrix updating step. They are described next.
\begin{figure*}
\centering
\includegraphics[height=4.5cm]{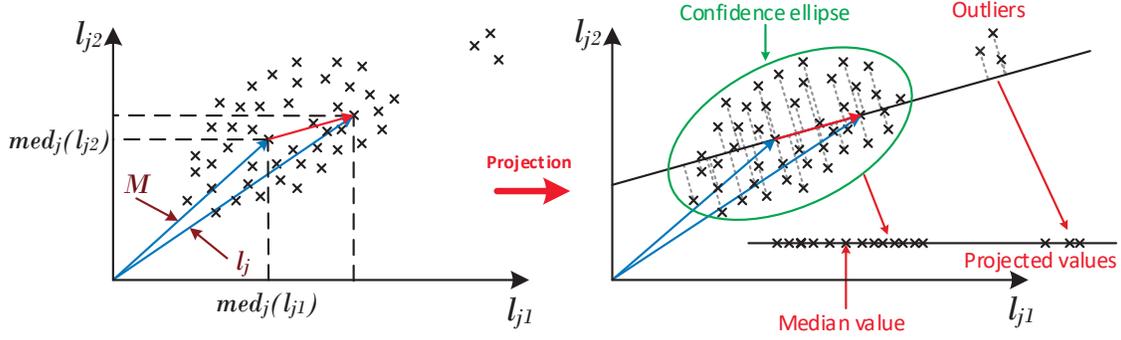}
\caption{Application of the projection statistics to the matrix $\bm{Z}_k$ for detecting outliers in a two-dimensional dataset that consists of the innovation vectors and the predicted state vectors.}
\label{Fig.PS}
\end{figure*}
\vspace{-0.3cm}
\subsection{Derivation of the Batch-Mode Regression Form}
In this subsection, we first show the equivalence of statistical linearization and the unscented transformation using sigma points. We then derive the proposed batch-mode regression form. The former claim is presented in the following Theorem:
\begin{theorem}
Given the state estimate vector ${\bm{\widehat x}_{k - 1\left| {k - 1} \right.}}$ and its associated covariance matrix ${\bm{P} _{k - 1\left| {k - 1} \right.}^{xx}}$, statistical linear regression applied to an arbitrary nonlinear function $\bm{g}(\bm{x})$ yields results that are equivalent to those of the unscented transformation using the sigma points generated according to (\ref{Eq:sigma_points1}).
\end{theorem}
\begin{proof}
Consider a nonlinear function $\bm{y}=\bm{g}(\bm{x})$ evaluated in 2$n$ points, i.e., $\left( {{\bm{\chi}_i},{\bm{\gamma}_i}} \right)$, where ${\bm{\gamma}_i} = \bm{g}\left( {{\bm{\chi}_i}} \right)$ for $i$= 1,..., 2$n$. Assuming that the nonlinear function is statistically linearized as $\bm{y}=\bm{Ax}+\bm{b}+\bm{\zeta}$, the objective is to find $\bm{\widehat A}$ and $\bm{\widehat b}$ so that the point-wise linearization error $\bm{\zeta}_i$ is minimized, i.e.,
\begin{equation}
\left\{ {\bm{\widehat A, \widehat b}} \right\} = \arg \min \sum\limits_{i = 1}^{2n} {{w_i}{\bm{\zeta}_i}^T{\bm{\zeta}_i}},
\label{Eq:statistical_linearization_objective_function}
\end{equation}
where $\bm{\zeta}_i=\bm{\gamma}_i-(\bm{A}\bm{\chi}_i+\bm{b})$. By taking the derivative of the objective function with respect to $\bm{A}$ and $\bm{b}$ and let them equal to zero, respectively, we obtain
\begin{equation}
\bm{b} = \bm{\overline y}  - \bm{\widehat A}\bm{\overline x},
\label{Eq:statistical_linearization_b}
\end{equation}
\begin{equation}
\bm{\widehat A}=\bm{P}_{xy}^T\bm{P}_{xx}^{-1},
\label{Eq:statistical_linearization_A}
\end{equation}
where $\bm{\overline x}  = \sum\limits_{i = 1}^{2n} {{w_i}{\bm{\chi _i}}}$; $\bm{\overline y} = \sum\limits_{i = 1}^{2n} {{w_i}\bm{g}\left( {{\bm{\chi _i}}} \right)}  = \sum\limits_{i = 1}^{2n} {{w_i}{\bm{\gamma}_i}}$; ${\bm{P}_{xx}} = \sum\limits_{i = 1}^{2n} {{w_i}\left( {{\bm{\chi}_i} - \bm{\overline x} } \right)} {\left( {{\bm{\chi}_i} - \bm{\overline x} } \right)^T}$; ${\bm{P}_{xy}} = \sum\limits_{i = 1}^{2n} {{w_i}\left( {{\bm{\chi}_i} - \bm{\overline x} } \right)} {\left( {{\bm{\gamma}_i} - \bm{\overline y} } \right)^T}$. Then, the estimation error covariance matrix is calculated as
\begin{equation}
\begin{array}{l}
{\bm{P}_{\zeta \zeta }} = \sum\limits_{i = 1}^{2n} {{w_i}{\bm{\widehat \zeta}_i}{\bm{\widehat\zeta}_i}^T} \\
{\rm{\quad   }} = \sum\limits_{i = 1}^{2n} {{w_i}\left( {{\bm{\gamma}_i} - \bm{\overline y}  - \bm{\widehat A}\left( {{\bm{\chi}_i} - \bm{\overline x} } \right)} \right){{\left( {{\bm{\gamma}_i} - \bm{\overline y}  - \bm{\widehat A}\left( {{\bm{\chi}_i} - \bm{\overline x} } \right)} \right)}^T}} \\
{\rm{ \quad  }} = {\bm{P}_{yy}} - \bm{\widehat A}{\bm{P}_{xx}}{\bm{\widehat A}^T} = {\bm{P}_{yy}} - \bm{P}_{xy}^T\bm{P}_{xx}^{ - 1}{\bm{P}_{xy}},
\end{array}
\label{Eq:statistical_linearization_covariance_matrix}
\end{equation}
where ${\bm{P}_{yy}} = \sum\limits_{i = 1}^{2n} {{w_i}\left( {{\bm{\gamma}_i} - \bm{\overline y}} \right)} {\left( {{\bm{\gamma}_i} - \bm{\overline y}} \right)^T}$.
Now, by taking the expectation and the outer product of the statistical linearized model, respectively, we obtain the posterior statistics given by
\begin{equation}
\bm{\widehat y}=\bm{\widehat A}\bm{\overline x}+\sum\limits_{i = 1}^{2n} {{w_i}{\bm{\chi} _i}}-\bm{\widehat A}\bm{\overline x}=\sum\limits_{i = 1}^{2n} {{w_i}{\bm{\chi} _i}},
\label{Eq:statistical_linearization_posterior_mean_statistics}
\end{equation}
\begin{equation}
\begin{array}{l}
{\bm{P}_{yy}} = \bm{\widehat A}{\bm{P}_{xx}}{\bm{\widehat A}^T} + {\bm{P}_{\zeta \zeta }}\\
{\rm{ \quad   }} = \bm{P}_{xy}^T\bm{P}_{xx}^{ - 1}{\bm{P}_{xy}} + \sum\limits_{i = 1}^{2n} {{w_i}\left( {{\bm{\gamma}_i} - \bm{\overline y} } \right)} {\left( {{\bm{\gamma}_i} - \bm{\overline y }} \right)^T} - \bm{P}_{xy}^T\bm{P}_{xx}^{ - 1}{\bm{P}_{xy}}\\
{\rm{ \quad   }} = \sum\limits_{i = 1}^{2n} {{w_i}\left( {{\bm{\gamma}_i} - \bm{\overline y} } \right)} {\left( {{\bm{\gamma}_i} - \bm{\overline y }} \right)^T},
\end{array}
\label{Eq:statistical_linearization_posterior_covariance_statistics}
\end{equation}
which are the same expressions as those obtained by applying the unscented transformation to the nonlinear function $\bm{y}=\bm{g}(\bm{x})$. Thus, the proof is completed.
\end{proof}
\begin{remark}
In statistical linearization, $\bm{\widehat A}$ is no longer the Jacobian matrix of $\bm{g}(\bm{x})$ at a given point. The error covariance matrix $\bm{P}_{\zeta\zeta}$ is used to compensate the linearization errors of the higher order Taylor series expansion terms. This is however explicitly contained in the unscented transformation process.
\end{remark}

By applying statistical linearization to the nonlinear system process model, we obtain the predicted state vector ${\bm{\widehat x}_{k\left| {k - 1} \right.}}$ along with its covariance matrix $\bm{P}_{k\left| {k - 1} \right.}^{xx}$. We define ${\bm{\widehat x}_{k\left| {k - 1} \right.}} = {\bm{x}_k} - {\bm{\delta} _k}$, where $\bm{x}_k$ is the true state vector; $\bm{\delta}_k$ is the prediction error; and $\mathbb{E}\left[ {{\bm{\delta } _k}\bm{\delta} _k^T} \right] = {\bm{P} _{k\left| {k - 1} \right.}^{xx}}$. Then, statistical linearization can be applied to the nonlinear observation equation, yielding
\vspace{-0.2cm}
\begin{equation}
{\bm{z}_k} = {\bm{H}_k}\left( {{\bm{x}_k} - {{\bm{\widehat x}}_{k\left| {k - 1} \right.}}} \right) + \bm{h}\left( {{{\bm{\widehat x}}_{k\left| {k - 1} \right.}}} \right) + {\bm{\nu}_k} + {\bm{\varepsilon}_k},
\label{Eq:statistical_linearization_observation_function}
\end{equation}
where ${\bm{H}_k} = {( {\bm{P}_{_{k\left| {k - 1} \right.}}^{xz}} )^T}(\bm{P}_{_{k\left| {k - 1} \right.}}^{xx})^{-1}$, which is no longer a Jacobian matrix. Here, the covariance of the statistical linearization error term is $\bm{\widetilde{R}}_k=\mathbb{E}\left[ {{\bm{\nu}_k}{\bm{\nu}_k^T}} \right] = \bm{P}_{_{k\left| {k - 1} \right.}}^{zz} - {( {\bm{P}_{_{k\left| {k - 1} \right.}}^{xz}} )^T}{\bm{P}_{k\left| {k - 1} \right.}^{xx}}\bm{P}_{_{k\left| {k - 1} \right.}}^{xz}$, where $\bm{P}_{_{k\left| {k - 1} \right.}}^{zz}$ and ${\bm{P}_{_{k\left| {k - 1} \right.}}^{xz}}$ are two covariance matrices that are calculated by following the same steps as those of the UKF. By processing the predictions and the observations simultaneously, we get the following batch-mode regression form:
\begin{equation}
\bigg[ {\begin{array}{*{10}{c}}
{{\bm{z}_k} + {\bm{H}_k}{{\bm{\widehat x}}_{k| {k - 1}}} - \bm{h}( {{{\bm{\widehat x}}_{k| {k - 1}}}})}\\
{{{\bm{\widehat x}}_{k|{k - 1}}}}
\end{array}}\bigg] = \bigg[ {\begin{array}{*{10}{c}}
{{\bm{H}_k}}\\
\bm{I}
\end{array}}\bigg]{\bm{x}_k} + \bigg[ {\begin{array}{*{10}{c}}
{{\bm{\nu}_k} + {\bm{\varepsilon}_k}}\\
{ - {\bm{\delta}_k}}
\end{array}}\bigg]
\label{Eq:batchmodel}
\end{equation}
which can be rewritten in a compact form as
\begin{equation}
{\bm{\widetilde z}_k} = {\bm{\widetilde H}_k}{\bm{x}_k} + {\bm{\widetilde e}_k},
\label{Eq:batchmodelcompact}
\end{equation}
and the error covariance matrix is given by
\begin{equation}
\bm{W}_k=\mathbb{E}\left[ {{{\bm{\widetilde e}}_k}\bm{\widetilde e}_k^T} \right] = \left[ {\begin{array}{*{20}{c}}
{{\bm{\Sigma} _{k\left| {k - 1} \right.}}}&\bm{0}\\
\bm{0}&{{\bm{P} _{k\left| {k - 1} \right.}^{xx}}}
\end{array}} \right] = {\bm{S}_k}\bm{S}_k^T,
\label{Eq:batchmodelcompactcovariance}
\end{equation}
where ${\bm{\Sigma}_{k| {k - 1}}} = \mathbb{E}[ {( {{\bm{\nu} _k} + {\bm{\varepsilon}_k}}){{( {{\bm{\nu}_k} + {\bm{\varepsilon}_k}})}^T}}]=\bm{R}_k+\bm{\widetilde{R}}_k$; $\bm{I}$ is an identity matrix; $\bm{S}_k$ is calculated by the Cholesky decomposition technique.
\begin{theorem}
The weighted least squares estimator of the batch-mode regression form (\ref{Eq:batchmodelcompact})
yields an estimated state vector ${{{\widehat x}_{k|{k}}}}$ and its associated covariance matrix $\bm{P}_{_{k|{k}}}^{xx}$ that are equivalent to those of the UKF.
\end{theorem}
\begin{proof}
It is well-known that the state estimate of (\ref{Eq:batchmodelcompact}) using the weighted least squares estimator is given by
\begin{equation}
{\bm{\widehat x}_{k\left| k \right.}} = {\left( {\bm{\widetilde H}_{_k}^T{\bm{W}_k}{\bm{\widetilde H}_k}} \right)^{ - 1}}\bm{\widetilde H}_{_k}^T{\bm{W}_k}{\bm{\widetilde z}_k},
\label{Eq:LSE_estimate}
\end{equation}
with the covariance matrix $\bm{P}_{_{k\left| {k} \right.}}^{xx}={\left( {\bm{\widetilde H}_{_k}^T{\bm{W}_k}{\bm{\widetilde H}_k}} \right)^{ - 1}}$. By applying an algebraic substitution and using the matrix inversion lemma, we get
\begin{equation}
\begin{array}{l}
\bm{P}_{_{k\left| k \right.}}^{xx} = {\left( {\bm{H}_k^T\bm{R}_k^{ - 1}{\bm{H}_k} + {{\left( {\bm{P}_{_{k\left| {k - 1} \right.}}^{xx}} \right)}^{ - 1}}} \right)^{ - 1}}\\
{\rm{}} = \bm{P}_{_{k\left| {k - 1} \right.}}^{xx} - \bm{P}_{_{k\left| {k - 1} \right.}}^{xx}\bm{H}_k^T{\left( {{\bm{H}_k}\bm{P}_{_{k\left| {k - 1} \right.}}^{xx}\bm{H}_k^T + {\bm{R}_k}} \right)^{ - 1}}{\bm{H}_k}\bm{P}_{_{k\left| {k - 1} \right.}}^{xx}\\
{\rm{}} = \left( {\bm{I} - {\bm{K}_k}{\bm{H}_k}} \right)\bm{P}_{_{k\left| {k - 1} \right.}}^{xx} = \bm{P}_{_{k\left| {k - 1} \right.}}^{xx} - {\bm{K}_k}\bm{P}_{_{k\left| {k - 1} \right.}}^{zz}\bm{K}_{_k}^T,
\end{array}
\label{Eq:LSE_estimate_covariance}
\end{equation}
where the gain matrix is expressed as
\begin{equation}
{\bm{K}_k}=\bm{P}_{_{k| {k - 1}}}^{xx}\bm{H}_k^T{( {{\bm{H}_k}\bm{P}_{_{k| {k - 1}}}^{xx}\bm{H}_k^T + {\bm{R}_k}})^{ - 1}}=\bm{P}_{_{k| {k - 1}}}^{xz}{( {\bm{P}_{_{k| {k - 1}}}^{zz}})^{ - 1}}
\label{Eq:LSE_estimate_gain_matrix}
\end{equation}
Thus, we can conclude that the estimation error covariance is identical to that of the UKF in (\ref{Eq:ukf_covariance_matrix_updating}). By applying similar substitutions and using the matrix inversion lemma, we can also show that the estimated state vector is given by
\begin{equation}
{\bm{\widehat x}_{k\left| k \right.}} = {\bm{\widehat x}_{k\left| {k - 1} \right.}} + {\bm{K}_k}\left( {{\bm{z}_{k}} - {{\bm{\widehat z}}_{k\left| {k - 1} \right.}}} \right),
\label{Eq:LSE_filtering}
\end{equation}
which completes the proof.
\end{proof}
\begin{remark}
In the literature, a few Huber estimator-based robust UKF methods have been proposed and applied to various applications in signal processing, target tracking, to name a few \cite{Chang2012,Chang2013,Chang2015}. However, in their developed regression models, $\bm{\nu}_k$ that compensates higher order Taylor series expansion error terms is neglected completely. As a consequence, the estimation results are biased. In addition, they are unable to handle innovation outliers and filter out non-Gaussian measurement noise.
\end{remark}
\vspace{-0.3cm}
\subsection{Robust Prewhitening}
Before carrying out a robust regression, we uncorrelate the state prediction errors of the batch-mode regression form. This is achieved by pre-multiplying $\bm{S}_k^{ - 1}$ on both sides of (\ref{Eq:batchmodelcompact}), yielding
\begin{equation}
\bm{S}_k^{ - 1}{\bm{\widetilde z}_k} = \bm{S}_k^{ - 1}{\bm{\widetilde H}_k}{\bm{x}_k} + \bm{S}_k^{ - 1}{\bm{\widetilde e}_k},
\label{Eq:prewhiteningbatchmodel}
\end{equation}
which can be further organized to the compact form
\begin{equation}
{\bm{y}_k} = {\bm{C}_k}{\bm{x}_k} + {\bm{\xi} _k},
\label{Eq:finalbatchmodel}
\end{equation}
where $\mathbb{E}[{\bm{\xi} _k}{\bm{\xi}_k}^T]=\bm{I}$.
\begin{figure*}
\centering
  \mbox{{\label{subfig:a} \includegraphics[width=8.2cm]{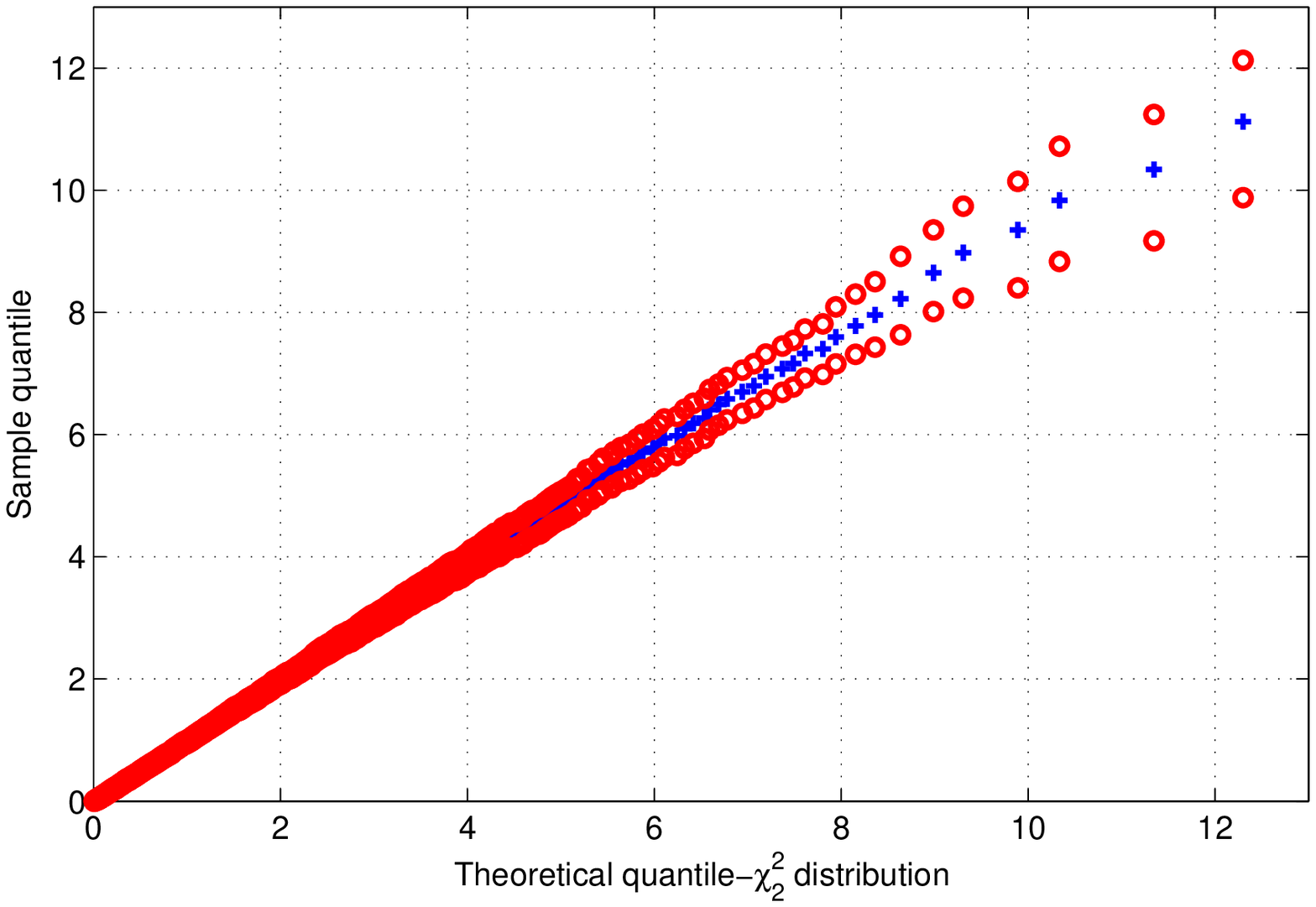}}}
  \mbox{{\label{subfig:b} \includegraphics[width=7cm]{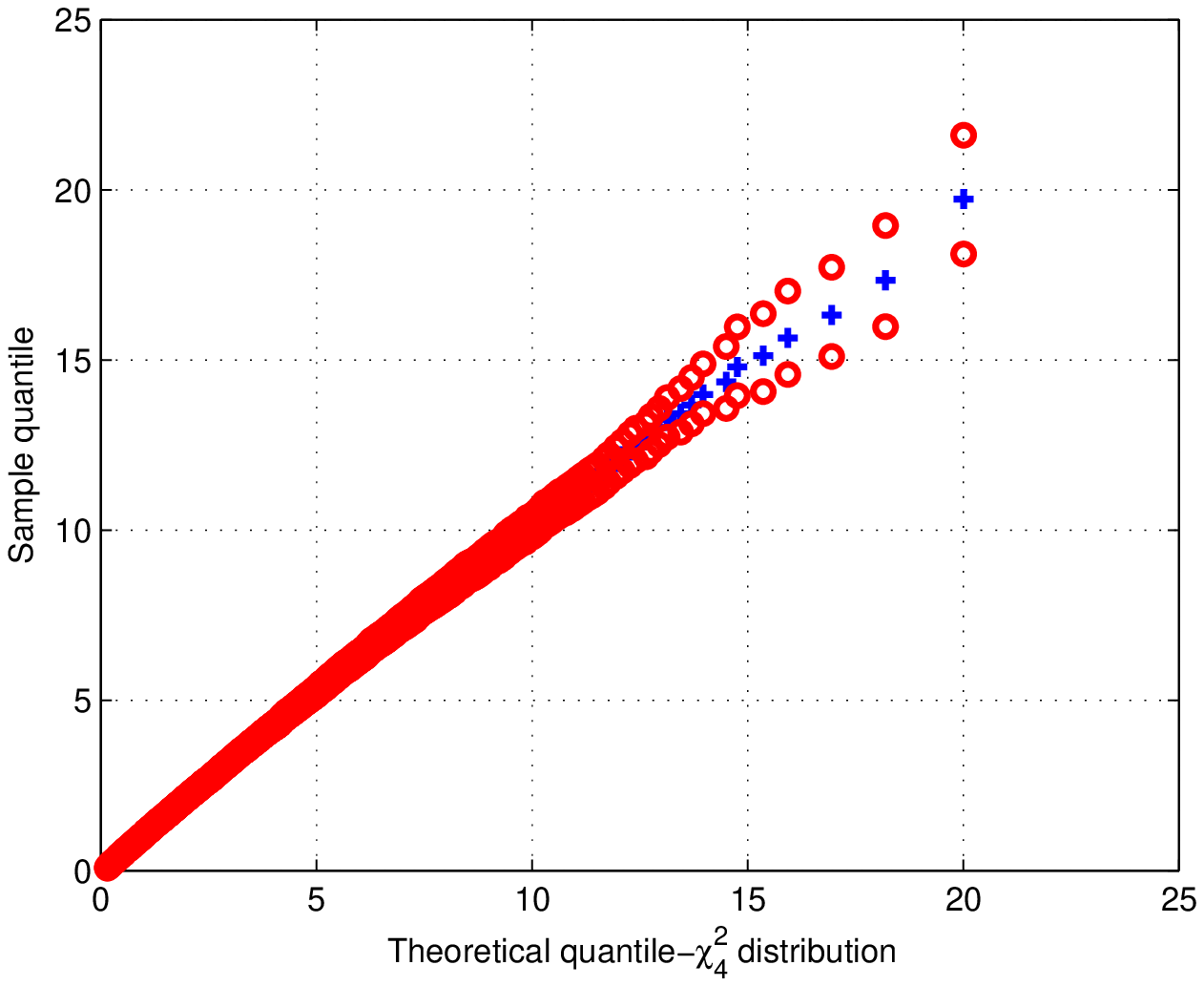}}}
\caption{Q-Q plots of the sample quantiles of the PS vs. the corresponding quantiles of the $\chi_2^2$ and $\chi_4^2$ distributions, where (a) and (b) represent Q-Q plots of PS with Gaussian and Laplace noise, respectively.}
\label{Fig.PS_distributions}
\end{figure*}
However, if outliers occur, the application of $\bm{S}_k^{ - 1}$ will corrupt the prewhitening \cite{Lmili2010}. To overcome this problem, we first detect the outliers and calculate the weights using the projection statistics (PS) \cite{Lmili2010,Lmili1996}. Those weights will be incorporated in the objective function that is defined in the proposed GM-estimator shown in Section III-C. Now, we describe the procedures used to calculate the weights. We apply the PS to a 2-dimensional matrix $\bm{Z}_k$ that contains serially correlated samples of the innovations and of the predicted state variables. Note that the innovation vector is defined as the difference between the observations and their associated predicted values at the previous step. Formally, we have
\vspace{-0.2cm}
\begin{equation}
\boldsymbol{Z}_k=
\left[ {\begin{array}{*{20}{c cc}}
\boldsymbol{z}_{k-1} - \boldsymbol{h}(\boldsymbol{\widehat{x}}_{k-1|k-2})   & \boldsymbol{z}_{k} - \boldsymbol{h}(\boldsymbol{\widehat{x}}_{k|k-1})\\
\boldsymbol{\widehat{x}}_{k-1|k-2}                                          & \boldsymbol{\widehat{x}}_{k|k-1}
\end{array}}
\right],
\label{Z_matrix}
\end{equation}
where $\boldsymbol{z}_{k-1} - \boldsymbol{h}(\boldsymbol{\widehat{x}}_{k-1|k-2})$ and $\boldsymbol{z}_{k} - \boldsymbol{h}(\boldsymbol{\widehat{x}}_{k|k-1})$ are the innovation vectors while $\widehat{x}_{k-1|k-2}$ and $\widehat{x}_{k|k-1}$ are the predicted state vectors at time instants $k$-1 and $k$, respectively. We may also apply the PS to higher dimensional samples, but we found that 2 dimensions are enough to identify outliers. The PS values of the predictions and of the innovations are separately calculated because the values taken by the former and the latter are centered around different points.

The PS of the $j$th row vector, $\bm{l}_{j}^T$, of the predictions (respectively the innovations) in $\bm{Z}_k$ is defined as the maximum of the standardized projections of all the $\bm{l}_{j}^T$'s on every direction $\bm{\ell}$ that originates from the coordinatewise medians of the predictions (respectively the innovations) and that passes through every data point, and where the standardized projections are based on the sample median and the median-absolute-deviation \cite{Lmili1996}. The implementation of the PS to detect outliers in matrix $\bm{Z}_k$ is displayed in Fig. \ref{Fig.PS}, while its mathematical expression is defined as \cite{Lmili1996}.
\vspace{-0.2cm}
\begin{equation}
P{S_j} = {\max _{\left\| \bm{\ell} \right\| = 1}}\frac{{\left| {\bm{l}_j^T\bm{\ell} - me{d_i}\left( {\bm{l}_i^T\bm{\ell} } \right)} \right|}}{{1.4826\;me{d_\kappa}\left| {\bm{l}_\kappa^T\bm{\ell}  - me{d_i}\left( {\bm{l}_i^T\bm{\ell}} \right)} \right|}},
\label{Eq:PS}
\end{equation}
where $i,j,\kappa=1,2,...,m+n$.

Once the PS values are calculated, they are compared to a statistical threshold to identify outliers. Extensive Monte Carlo simulations and Q-Q plots reveal that the probability distributions of the PS applied to $\bm{Z}_k$, whose data points obey bivariate Gaussian and Laplace probability distributions, follow chi-square distributions with degree of freedom 2 and 4, respectively (See Fig. \ref{Fig.PS_distributions}). This investigation allows us to apply statistical tests to the PS and to flag all the data points that satisfy $\text{PS}_i>\eta$ as outliers. The latter are downweighted via
\begin{equation}
{\varpi _i} = \min \left( {1,{\rm{ }}{{{d^2}} \mathord{\left/
 {\vphantom {{{d^2}} {PS_i^2}}} \right.
 \kern-\nulldelimiterspace} {PS_i^2}}} \right),
\label{Eq:PSdownweightfunction}
\end{equation}
where the parameter $d$ is set equal to 1.5 to yield good statistical efficiency at different distributions without increasing too much the bias induced by outliers. As an example, when the noise is assumed to be Laplacian, the PS obeys a chi-square distribution with 4 degrees of freedom. In that case, we can choose the statistical detection threshold $\eta$ as $\chi_{4,0.975}^2$ at a significance level of 97.5\%.

\vspace{-0.3cm}
\subsection{Robust Filtering and Solution}
To suppress the outliers and filter out thick-tailed non-Gaussian measurement noise, we develop a robust GM-estimator that minimizes the following objective function:
\begin{equation}
J\left( \bm{x}_k \right) = \sum\limits_{i = 1}^{m+n} {\varpi _i^2} \rho \left( {{r_{{S_i}}}} \right),
\label{Eq:objectivefunction}
\end{equation}
where $\varpi _i$ is calculated by (\ref{Eq:PSdownweightfunction}); ${r_{{S_i}}} = {{{r_i}} \mathord{\left/{\vphantom {{{\bm{r}_i}} {s{\varpi _i}}}} \right.\kern-\nulldelimiterspace} {s{\varpi _i}}}$ is the standardized residual; ${r_i} = {y_i} -\bm{c}_i^T \bm{\widehat x}$ is the residual, where $\bm{c}_i^T$ is the $i$th row vector of the matrix $\bm{C}_k$; $s = 1.4826\cdot b_m{\rm{\cdot}}\text{median}_i\left| {{r_i}}\right|$ is the robust scale estimate; $b_m$ is a correction factor to achieve unbiasedness for a finite sample of size $m+n$ at a given probability distribution; $\rho(\cdot)$ is the nonlinear function of ${r_{{S_i}}}$. In this paper, the convex Huber-$\rho$ function \cite{Huber1981} is adopted, that is
\begin{equation}
\rho \left( {{r_{{S_i}}}} \right) = \left\{ {\begin{array}{*{20}{c}}
{\frac{1}{2}r_{_{{S_i}}}^2,{\rm{ \quad\quad\quad\quad  \quad            for }}\left| {{r_{{S_i}}}} \right| < \lambda}\\
{\lambda\left| {{r_{{S_i}}}} \right| - {{{\lambda^2}} \mathord{\left/
 {\vphantom {{{\lambda^2}} 2}} \right.
 \kern-\nulldelimiterspace} 2},{\rm{ \quad\; }}elsewhere}
\end{array}} \right.,
\label{Eq:huberfunction}
\end{equation}
where the parameter $\lambda$ between the quadratic and the linear segment of $\rho(\cdot)$ is typically chosen between 1.5 to 3 in the literature.

To minimize (\ref{Eq:objectivefunction}), one takes its partial derivative with respect to $\bm{x}_k$ and sets it equal to zero, yielding
\begin{equation}
\frac{{\partial J\left( \bm{x}_k \right)}}{{\partial \bm{x}_k}} = \sum\limits_{i = 1}^{m+n} { - \frac{{{\varpi _i}{\bm{c}_i}}}{s}\psi \left( {{r_{{S_i}}}} \right)}  = \bm{0},
\label{Eq:objectivefunctionpartialderivative}
\end{equation}
where $\psi \left( {{r_{{S_i}}}} \right) = {{\partial \rho \left( {{r_{{S_i}}}} \right)} \mathord{\left/
 {\vphantom {{\partial \rho \left( {{r_{{S_i}}}} \right)} {{ r_{{S_i}}}}}} \right.
 \kern-\nulldelimiterspace} \partial{{r_{{S_i}}}}}$is the so-called $\psi$-function. By dividing and multiplying the standardized residual $r_{{S_i}}$ to both sides of (\ref{Eq:objectivefunctionpartialderivative}) and putting it in a matrix form, we get
\begin{equation}
{\bm{C}_k^T}\bm{\widehat{Q}}\left( {\bm{y}_k -\bm{C}_k \bm{x}_k } \right) = \bm{0},
\label{Eq:objectivefunctionpartialderivativematrixform}
\end{equation}
where $\bm{\widehat Q} =$diag${\left( {q\left( {{r_{{S_i}}}} \right)} \right)}$  and $q\left( {{r_{{S_i}}}} \right) = {{\psi \left( {{r_{{S_i}}}} \right)} \mathord{\left/{\vphantom {{\psi \left( {{r_{{S_i}}}} \right)} {{r_{{S_i}}}}}} \right.\kern-\nulldelimiterspace}{{r_{{S_i}}}}}$.

By using the IRLS algorithm \cite{Lmili2007,Hampel1986}, the state vector correction at the $j$ iteration is calculated through
\begin{equation}
\Delta \bm{\widehat x}_{k\left| k \right.}^{\left( {j + 1} \right)} = {\left( {\bm{C}_k^T{\bm{\widehat Q}^{\left( j \right)}}{\bm{C}_k}} \right)^{ - 1}}\bm{C}_k^T{\bm{\widehat Q}^{\left( j \right)}}\bm{y}_k,
\label{Eq:IRLS}
\end{equation}
where $\Delta \bm{\widehat x}_{k\left| k \right.}^{\left( {j + 1} \right)}=\bm{\widehat x}_{k\left| k \right.}^{\left( {j + 1} \right)}-\bm{\widehat x}_{k\left| k \right.}^{\left( {j} \right)}$. The algorithm converges when ${\left\| {\Delta \bm{\widehat x}_{k\left| k \right.}^{\left( {j + 1} \right)}} \right\|_\infty } \le {10^{ - 2}}$.

\subsection{Asymptotic Error Covariance Matrix of the GM-UKF State Estimates}
Upon convergence of the iterative algorithm, the error covariance matrix $\bm{P}_{k\left|{k}\right.}^{xx}$ is updated so that the state prediction for the next step can be performed. To this end, consider the $\epsilon$-contamination model $G = \left( {1 - \epsilon } \right)\Phi + \epsilon {\Delta _r}$, where $G$ and $\Phi$ are the contaminated and the true cumulative probability distribution function of the residual, respectively; ${\Delta _r}$ is the point mass to model outliers or unknown non-Gaussian distributions. The error covariance matrix is updated based on the following theorem:
\vspace{-0.2cm}
\begin{theorem}
Let $\bm{T}(\cdot)$ be the functional form of the GM-estimator with a bounded $\psi(\cdot)$ function and $\Phi_\alpha$ be the empirical cumulative probability distribution function, then
\begin{equation}
\sqrt \alpha ( {\bm{T}( {{\Phi_\alpha }}) - \bm{T}(\Phi)})\mathop  \to \limits^d \mathcal{N}( {\bm{0},\bm{P}_{k| k }^{xx}}),
\label{Eq:Gaussian_convergence}
\end{equation}
where $\alpha=m+n$; $\mathop \to \limits^d$ means convergence in probability distribution; $\bm{P}_{k| k }^{xx}=\mathbb{E}[ {\bm{IF}(\bm{x};\Phi,\bm{T}) \cdot \bm{I{F}}(\bm{x};\Phi,\bm{T})^T}]$ with the influence function $\bm{IF}(\bm{x};\Phi,\bm{T})$ evaluated at $\Phi$.
\end{theorem}
\begin{proof}
By taking a first-order Taylor series expansion of the functional form of the estimator $\bm{T}$ with respective to $\Phi$, we get
\begin{equation}
\bm{T}( {{\Phi_\alpha }}) = \bm{T}(\Phi) + {\bm{T}^{\prime}}( {{\Phi_\alpha } - \Phi}) + {\mathop{\rm Rem}\nolimits}( {{\Phi_\alpha } -\Phi}),
\label{Eq:Taylor_expansion}
\end{equation}
which can be reorganized into the following form by multiplying $\sqrt \alpha$ on both sides of the equality:
\begin{IEEEeqnarray}{lll}
\sqrt \alpha  \left( {\bm{T}\left( {{\Phi_\alpha }} \right) - \bm{T}\left(\Phi\right)} \right)\nonumber\\
 = \sqrt \alpha  {\bm{T}^{\prime}}\left( {{\Phi_\alpha } -\Phi} \right) + \sqrt \alpha  {\mathop{\rm Rem}\nolimits}\left( {{\Phi_\alpha } -\Phi} \right)\\
 = \sqrt \alpha\int {\bm{IF}({\bm{x};\Phi,\bm{T}})} d({{\Phi_\alpha }-\Phi}) + \sqrt \alpha  {\mathop{\rm Rem}\nolimits}( {{\Phi_\alpha  -\Phi}})\\
  = \sqrt \alpha\int {\bm{IF}({\bm{x};\Phi,\bm{T}})} d{{\Phi_\alpha}} + \sqrt \alpha  {\mathop{\rm Rem}\nolimits}( {{\Phi_\alpha  -\Phi}})\\
 = \frac{1}{\sqrt \alpha}\sum\limits_{i = 1}^\alpha  {\bm{IF}\left({{x_i};\Phi,\bm{T}} \right)}  + \sqrt \alpha  {\mathop{\rm Rem}\nolimits}\left( {{\Phi_\alpha } -\Phi} \right),
\label{Eq:Taylor_functional}
\end{IEEEeqnarray}
where the definition of the influence function is applied to yield (38) to (39); by virtue of Fisher consistency at the distribution $\Phi$, that is, $\int {\bm{IF}({\bm{x};\Phi,\bm{T}})} d{\Phi}=\bm{0}$, (39) reduces to (40); finally, by using the property of the empirical cumulative probability distribution function, we have
\begin{equation}
\int {\bm{IF}({\bm{x};\Phi,\bm{T}})} d{{\Phi_\alpha}}=\frac{1}{\alpha}\sum\limits_{i = 1}^\alpha{\bm{IF}\left({{x_i};\Phi,\bm{T}} \right)},
\label{Eq:empirical_distribution}
\end{equation}
yielding (40) to (41).

Following the work of Fernholz \cite{Fernholz1983}, we can show that
\begin{equation}
\sqrt \alpha  {\mathop{\rm Re}\nolimits} m\left( {{\Phi_\alpha }-\Phi} \right)\mathop  \to \limits^p 0,
\label{Eq:Rem_convergence}
\end{equation}
where $\mathop  \to \limits^p$ means probability convergence. Therefore, by applying the central limit theorem and Slutsky's lemma to (\ref{Eq:Taylor_functional}), it follows that
\begin{equation}
\sqrt \alpha ( {\bm{T}( {{\Phi_\alpha }}) - \bm{T}(\Phi)})\mathop  \to \limits^d \mathcal{N}( {\bm{0},\bm{P}_{k| k }^{xx}}),
\label{Eq:Gaussian_convergence}
\end{equation}
where $\bm{P}_{k| k }^{xx}=\mathbb{E}[ {\bm{IF}(\bm{x};\Phi,\bm{T}) \cdot \bm{I{F}}(\bm{x};\Phi,\bm{T})^T}]$.
\end{proof}

\emph{Discussion}: The UKF is able to provide good results only when the process and observation noises obey a Gaussian distribution \cite{Julier2000}. Indeed, in that case the filtered state vector ${\bm{\widehat {x}}_{k| {k}}}$ is Gaussian and the mean and covariance matrix of ${\bm{\widehat {x}}_{k| {k}}}$ can be accurately estimated by the sample mean and the sample covariance matrix of the sigma points. However, this property no longer holds true if the Gaussianity assumption of the noises is violated. In that case, the state estimate vector ${\bm{\widehat {x}}_{k| {k}}}$ obtained from the UKF is significantly biased due the filter lack of statistical robustness to thick-tailed non-Gaussian noise. By contrast, our GM-UKF guarantees the asymptotic Gaussianity of ${\bm{\widehat {x}}_{k| {k}}}$ for thick-tailed non-Gaussian noises and yields reliable state estimates with good statistical efficiency.

\begin{corollary}
Assume that the system process noise is contaminated about a Gaussian distribution. Then, the data points defined by the row vectors of the matrix $\bm{Z}_k$ follow asymptotically a Gaussian distribution.
\end{corollary}
\begin{proof}
From the definition of the matrix $\bm{Z}_k$ given by (\ref{Z_matrix}) and Theorem 3, we can see that the predicted state vector is roughly Gaussian. Furthermore, under the assumption that the minority of the measurements obey a thick-tailed non-Gaussian distribution, the innovation vectors can be shown to be approximately Gaussian. From this, we conclude that $\bm{Z}_k$ is asymptotically Gaussian.
\end{proof}

Let's now derive the $\bm{IF}(\bm{x};\Phi,\bm{T})$ of our GM-UKF at the cumulative probability distribution $\Phi$.

\begin{corollary}
The total influence function of the GM-UKF defined by (\ref{Eq:objectivefunctionpartialderivative}) using the regression model (\ref{Eq:finalbatchmodel}) is expressed as
\begin{equation}
\bm{IF}(\bm{x};\Phi,\bm{T}) = {\left[ {\int {\frac{1}{s}{\psi ^{'}}\left( {{r_{{S_i}}}} \right)\bm{C{C}^T}\left| {_{_{_{T\left(\Phi \right)}}}d\Phi} \right.} } \right]^{ - 1}}\varpi \bm{C}\psi \left( {{r_{{S_i}}}} \right).
\label{Eq:IFfinal}
\end{equation}
\end{corollary}

\begin{proof}
In our previous work \cite{Junbo_GMIEKF2016}, the total influence function of a GM-estimator based on a nonlinear regression model given by $\bm{y}=\bm{\varphi}(\bm{x})+\bm{\xi}$ is expressed as
\begin{equation}
\normalsize
\begin{array}{l}
\bm{IF}(\bm{x};\Phi,\bm{T})\\
 = {\left( {\int {\left\{ {\frac{{{\psi ^{'}}\left( {{\bm{r}_{{S_i}}}} \right)}}{s}\frac{{\partial \bm{\varphi \left( {x} } \right)}}{{\partial \bm{x}}}{{\frac{{\partial \bm{\varphi \left( {x}} \right)}}{{\partial \bm{x}}}}^T} - \varpi \psi \left( {{\bm{r}_{{S_i}}}} \right)\bm{D}} \right\}\left| {_{T\left( \Phi \right)}}d\Phi \right.} } \right)^{ - 1}}\\
{\rm{ \quad }} \cdot \varpi \frac{{\partial \bm{\varphi \left( {x}} \right)}}{{\partial \bm{x}}}\psi \left( {{\bm{r}_{{S_i}}}} \right),
\end{array}
\label{Eq:Genral_IF}
\end{equation}
where $\bm{D}={\frac{{{\partial ^2}\bm{\varphi (x)}}}{{\partial {x_i}\partial {x_j}}}}$ is the Hessian matrix of $\bm{\varphi (x)}$.
Since we have $\bm{\varphi}(\bm{x})=\bm{C}\bm{x}$ for the GM-UKF, (\ref{Eq:Genral_IF}) reduces to
\begin{equation}
\bm{IF}(\bm{x};\Phi,\bm{T}) = {\left[ {\int {\frac{1}{s}{\psi ^{'}}\left( {{r_{{S_i}}}} \right)\bm{C{C}^T}\left| {_{_{_{T\left( \Phi \right)}}}d\Phi} \right.} } \right]^{ - 1}}\varpi \bm{C}\psi \left( {{r_{{S_i}}}} \right).
\label{Eq:IFfinal}
\end{equation}
\end{proof}

Now, we are in a position to derive the covariance matrix $\bm{P}_{k| k }^{xx}$ from (\ref{Eq:Gaussian_convergence}). First, let us prove the following theorem:
\begin{theorem}
The sample variance of the robust scale estimator $s$ of the GM-standardized residuals tends to one as the number of observation tends to infinity.
\end{theorem}
\begin{proof}
By the law of large numbers, the distribution of the residuals tends to the Gaussian distribution, i.e., $\Phi\sim\mathcal{N}(\mu,\sigma^2)$. Since the median absolute deviation (MAD=$\frac{s}{{1.4826 \cdot {b_m}}}$) is a consistent estimator for the standard deviation $\sigma$ of a Gaussian distribution, we get
\begin{equation}
\begin{array}{l}
\frac{1}{2}{\rm{ = }}P\left( {\left| {X - \mu } \right| \le \frac{s}{{1.4826 \cdot {b_m}}}} \right) = P\left( {\frac{{\left| {X - \mu } \right|}}{\sigma } \le \frac{s}{{1.4826 \cdot {b_m} \cdot \sigma }}} \right)\\
 \;\;= 2\Phi \left( {\frac{s}{{1.4826 \cdot {b_m} \cdot \sigma }}} \right) - 1.
\end{array}
\end{equation}
Therefore, we obtain $s/\sigma = 1.4826 \cdot {b_m} \cdot {\Phi ^{ - 1}}\left( {\frac{3}{4}} \right) = {b_m} \to 1$ as $m$ tends to infinity, where $\Phi$ is the cumulative probability function of the standard Gaussian distribution. On the other hand, from the equation (\ref{Eq:finalbatchmodel}) and the fact that $\mathbb{E}[{\bm{\xi} _k}{\bm{\xi}_k}^T]=\bm{I}$, the residuals can be shown to actually follow the standard Gaussian distribution. Therefore, $\mathbb{E}_F[s^2]=s^2\to 1$.
\end{proof}

Finally, by Theorems 3 and 4, the asymptotic error covariance matrix of our GM-UKF at time sample $k$ is updated by
\begin{equation}
\begin{array}{l}
{\bm{P} _{k\left| k \right.}^{xx}}=\mathbb{E}[ {\bm{IF}(\bm{x};\Phi,\bm{T}) \cdot \bm{I{F}}(\bm{x};\Phi,\bm{T})^T}]\\
\quad\quad\; = \frac{{{\mathbb{E}_\Phi}\left[ {{\psi ^2}\left( {{r_{{S_i}}}} \right)} \right]}}{{{{\left\{ {{\mathbb{E}_\Phi}\left[ {{\psi ^{\prime}}\left( {{r_{{S_i}}}} \right)} \right]} \right\}}^2}}}{\left( {{\bm{C}_k^T}\bm{C}_k} \right)^{ - 1}}{\left( {{\bm{C}_k^T}{\bm{Q}_\varpi }\bm{C}_k} \right)^{}}{\left( {{\bm{C}_k^T}\bm{C}_k} \right)^{ - 1}}
\end{array}
\label{Eq:updatingSigmafinal}
\end{equation}
where ${\bm{Q}_\varpi} = diag\left( {\varpi _i^2} \right)$.

\subsection{Discussions on the Statistical Efficiency of the GM-UKF}
In this section, we discuss the statistical efficiency of our proposed GM-UKF under various probability distributions of the noise. Firstly, under Gaussian measurement noise, the outliers detected by the PS will be downweighted by the linear segment of the $\rho$-function while all the good measurements will be assigned weights equal or close to one since most of them will be processed by the quadratic segment of the $\rho$-function. As a result, the state estimator exhibits a high statistical efficiency. Secondly, if the measurement noise obeys a Laplace distribution, those measurements associated with the thick tails of that distribution will have standardized residuals corresponding to the linear segment of the $\rho$-function. This means that for them, the GM-estimator behaves like the least absolute value estimator; since the latter is the maximum-likelihood estimator at that distribution, it will have a high asymptotic statistical efficiency. On the other hand, when the estimation error covariance matrix is updated, all the outliers with respect to the Gaussian distribution, which include the measurements associated with the tails of the Laplacian or the Cauchy distribution, will be heavily downweighted through the matrix ${\bm{Q}_\varpi}$, yielding bounded biases and variances in the state estimates.

\section{Conclusion}
In this first part of a two-part series paper, we present the fundamental theory of the proposed GM-UKF. We show first that the UKF estimates the state vector via a weighted least squares estimator under the Gaussianity assumption of the system process or measurement noises; consequently, it yields strongly biased state estimates when the noises follow non-Gaussian probability distributions, which is precisely the case when processing PMU measurements. By contrast, the state estimates and residuals of our GM-UKF are proved to be asymptotically Gaussian, allowing the sigma points to reliably approximate the mean and the covariance matrices of the predicted and corrected state vectors. Furthermore, by relying on the projection statistics and the GM-estimator, the proposed GM-UKF is able to suppress observation and innovation outliers while exhibiting high statistical efficiency of the state estimates. In addition, we derive the expression of the asymptotic error covariance matrix of the GM-UKF state estimates from the total influence function of the GM-estimator. In the companion paper, we will discuss the implementation of our GM-UKF in power systems and analyze its performance by carrying out extensive simulations under various scenarios.

\ifCLASSOPTIONcaptionsoff
  \newpage
\fi
\end{document}